\documentclass[a4paper,reqno]{amsart}
\usepackage{dsfont}
\usepackage{newcent} 
\usepackage[T1]{fontenc}
\usepackage{amsmath,amssymb,amsfonts,amsthm}
\usepackage{mathtools}
\usepackage{newlfont}
\usepackage{fancyhdr,fancyvrb}
\usepackage{graphicx}
\usepackage{nextpage}
\usepackage[dvipsnames,usenames]{color}
\usepackage[colorlinks=true, linkcolor=blue, citecolor=red]{hyperref}
\usepackage{tikz}
\usepackage{tikz-3dplot}
\usepackage{mathrsfs}
\usepackage{stmaryrd}

 \makeatletter
 \@namedef{subjclassname@2020}{%
 \textup{2020} Mathematics Subject Classification}
 \makeatother

\newcommand{\M}[1]{\mathbb{M}^{#1}}
\newcommand{\N}{\mathbb{N}}
\newcommand{\R}[1]{\mathbb{R}^{#1}}
\renewcommand{\S}[1]{\mathbb{S}^{#1}}

\newcommand{\bP}{\mathbf{P}}

\newcommand{\cA}{\mathcal A}

\newcommand{\cC}{\mathcal C}
\newcommand{\cD}{\mathcal D}
\newcommand{\cE}{\mathcal E}

\newcommand{\cH}{\mathcal H}

\newcommand{\cV}{\mathcal V}
\newcommand{\cW}{\mathcal W}
\newcommand{\cX}{\mathcal X}

\newcommand\bigw{\scalebox{.95}[1]{$\bigwedge$}}

\newcommand{\de}{\mathrm d}

\newcommand{\spann}{\operatorname{span}}

\newcommand{\eps}{\varepsilon}

\newcommand{\cof}{\operatorname{cof}}
\newcommand{\curv}{\operatorname{curv}}

\newcommand{\diam}{\operatorname{diam}}

\newcommand{\Id}{\operatorname{I}}

\newcommand{\spt}{\operatorname{spt}}
\newcommand{\tr}{\operatorname{tr}}

\newcommand{\e}{\varepsilon}
\usepackage{soul}

\newcommand{\llb}{\llbracket}
\newcommand{\rrb}{\rrbracket}

\newcommand{\red}[1]{\textcolor[rgb]{1,0,0}{#1}}

\newcommand{\blu}[1]{\textcolor[rgb]{0,0,0}{#1}}

\renewcommand{\div}{\operatorname{div}}
\renewcommand{\geq}{\geqslant}
\renewcommand{\leq}{\leqslant}
\newcommand{\wto}{\rightharpoonup}

\newcommand{\average}{{\mathchoice {\kern1ex\vcenter{\hrule
height.4pt width 8pt depth0pt}
\kern-11pt} {\kern1ex\vcenter{\hrule height.4pt width 4.3pt
depth0pt} \kern-7pt} {} {} }}

\newcommand{\res}{\mathop{\hbox{\vrule height 7pt width .5pt depth
0pt\vrule height .5pt width 6pt depth0pt}}\nolimits}

\mathchardef\emptyset="001F

\numberwithin{equation}{section}

\newtheorem{defin}{Definition}[section]
\newtheorem{remark}[defin]{Remark}
\newtheorem{theorem}[defin]{Theorem}
\newtheorem{lemma}[defin]{Lemma}
\newtheorem{proposition}[defin]{Proposition}

\title{Direct minimization of the Canham--Helfrich energy on generalized Gauss graphs}%
\author{Anna Kubin}
\address[A.~Kubin]{Dipartimento di Scienze Matematiche ``G.~L.~Lagrange'', Politecnico di Torino, Corso Duca degli Abruzzi, 24, 10129 Torino, Italy \&
Institute of Analysis and Scientific Computing, Technische Universit\"at Wien, Wiedner Haupstrasse 8-10, 1040 Vienna, Austria}
\email{anna.kubin@polito.it \& anna.kubin@asc.tuwien.ac.at}
\author{Luca Lussardi}
\address[L.~Lussardi]{Dipartimento di Scienze Matematiche ``G.~L.~Lagrange'', Politecnico di Torino, Corso Duca degli Abruzzi, 24, 10129 Torino, Italy}
\email{luca.lussardi@polito.it}
\author{Marco Morandotti}
\address[M.~Morandotti]{Dipartimento di Scienze Matematiche ``G.~L.~Lagrange'', Politecnico di Torino, Corso Duca degli Abruzzi, 24, 10129 Torino, Italy}
\email{marco.morandotti@polito.it}

\date{\today}


\makeindex

\begin{document}

\begin{abstract}
The existence of minimizers of the Canham--Helfrich functional in the setting of generalized Gauss graphs is proved. 
As a first step, the Canham--Helfrich functional, usually defined on regular surfaces, is extended to generalized Gauss graphs, then lower semicontinuity and compactness are proved under a suitable condition on the bending constants ensuring coerciveness; the minimization follows by the direct methods of the Calculus of Variations. Remarks on the regularity of minimizers and on the behavior of the functional in case there is lack of coerciveness are presented.
\end{abstract}

\maketitle

\tableofcontents



\allowdisplaybreaks

\section{Introduction} \label{sec:intro}
The mathematical modeling of biological membranes is an active field of research that has received much attention in the last half century starting with the pioneering works of Canham \cite{Canham1970} and Helfrich \cite{Helfrich1973}. 
They modeled the membranes as regular surfaces in the space and associate the equilibrium configurations with the minimum of an energy functional depending on the curvatures. 
If we denote by $M\subseteq \R{3}$ a two-dimensional, compact, and oriented {sub}manifold (with an understood choice of the normal vector $\nu\colon M\to\mathbb{S}^2$), by $H$ and $K$ its mean and Gaussian curvatures, respectively, and by $H_0$ a constant \emph{spontaneous curvature}, the so-called \emph{Canham--Helfrich energy functional} reads
\begin{equation}\label{100}
\cE(M)\coloneqq \int_M \big(\alpha_H(H(p)-H_0)^2-\alpha_K K(p)\big)\,\de\cH^{2}(p),
\end{equation}
where $\cH^{2}$ is the $2$-dimensional Hausdorff measure and 
$\alpha_H,\alpha_K>0$ 
are the \emph{bending constants}. 
These are physical, model-specific, constants and the range of possible values that they can take reveals to be crucial to determine the coercivity of the functional. We point out that 
the positivity of the constants and the minus sign between the two terms in the energy give competition between the two curvatures in \eqref{100}.

In the smooth case where $M$ is at least of class $\cC^2$, the curvatures $H$ and $K$ are given by the usual formulae
\begin{equation*}
H=\kappa_1+\kappa_2\quad\text{and}\quad K=\kappa_1\kappa_2,
\end{equation*}
$\kappa_i$ being the principal curvatures, with respect to which, \blu{when $H_0=0$}, the functional $\cE$ is homogeneous of degree two.

If $M$ is without boundary, one can invoke the Gauss--Bonnet theorem and obtain that the term involving $K$ gives a constant contribution (determined by the Euler characteristic $\chi(M)$ of $M$) to the energy, so that it can be neglected in view of the minimization of~$\cE$ among all surfaces with prescribed topology.
In this case, and under the further constraint that the spontaneous curvature vanishes, the functional $\cE$ reduces to the well known \emph{Willmore energy functional} \cite{KS2012,Riviere2008,Schygulla2012,Simon,Willmore1965} 
\begin{equation}\label{100W}
\cW(M)\coloneqq \int_M H^2(p)\,\de\cH^{2}(p).
\end{equation}

Both functionals $\cE$ and $\cW$ are geometric in nature, since they depend on 
geometric features of the surface $M$, and can be studied in a number of different contexts, according to the regularity requests on $M$.
Sobolev-type approaches to the minimization either of the Willmore functional (see \cite{KS2012} and the references therein, see also \cite{KMG2023,KMR2014,MS2023,RS2023}) or of the Canham--Helfrich functional (see, \emph{e.g.}, \cite{CMV,CV,H2013,H2015,LR,MS,W}) assume that $M$ has fixed topology, or even symmetry constraints. 
Aiming at considering more general surfaces, a successful approach is the one through varifolds \cite{LCY,Hutchinson}, see \cite{BLS,Eichmann2019,Eichmann2020}.
We point out that other frameworks are available in the study of geometric functionals: for instance, currents \cite{FF1960} have been used to tackle the minimization of the area functional. Despite not being suitable for the formulation of problems involving curvatures, due to their lack of an intrinsic notion of curvature, special classes of currents have been introduced to overcome this issue.
Nonetheless, it is possible to apply the technical tools of the theory of currents to the class of the so-called \emph{generalized Gauss graphs} \cite{AST}, which are motivated by a generalization of the graph of the Gauss map on smooth surfaces $M$. 
Instead of generalizing $M$ itself, this approach has the remarkable advantage to allow one to exploit the fact that the curvatures of $M$ are coded in its Gauss map, see Section~\ref{sec:GGG} for details.

We point out that the typical form in which the functional $\cE$ in \eqref{100} is found in the literature is 
\begin{equation}\label{100tilde}
E(M)\coloneqq\int_M \Big(\frac{a_H}2(H(p)-H_0)^2+a_K K(p)\Big)\,\de\cH^{2}(p),
\end{equation}
under the condition that $a_H>0$ and $a_K<0$ to ensure the competition between the two curvatures. In this context, it is required that 
\begin{equation}\label{costantitilde}
a_H>0\qquad\text{and}\qquad \frac{a_K}{a_H}\in(-2,0)
\end{equation}
in order to ensure both the coercivity and the lower semicontinuity of the functional \eqref{100tilde}; this condition is the same assumed in \cite[Theorem~1]{CV} and \cite[formula~(1.9)]{CMV} in the Sobolev setting, see also \cite{H2013,H2015}, whereas the more restrictive condition $-6a_H<5a_K<0$ is considered in \cite{BLS} in the varifold setting. 
We note that the typical physical range of the parameters is $-a_H\leq a_K\leq0$, see, \emph{e.g.}, \cite{BAUMGART20051067,BM2010,TKS1998}\footnote{Note that \cite{BAUMGART20051067} cites \cite{TKS1998}, but inverting numerator and denominator, by mistake.}, the case $a_K=0$ essentially reducing to the Willmore functional $\cW$ of \eqref{100W}.
Given the expression of the Canham--Helfrich functional $\cE$ in \eqref{100}, condition \eqref{costantitilde} reads
\begin{equation}\label{costanti}
4\alpha_H>\alpha_K>0.
\end{equation}

In this paper, we provide a suitable formulation of the Canham--Helfrich functional $\cE$ introduced in \eqref{100} in the class of generalized Gauss graphs and study three minimization problems. 
Our main results are Theorems~\ref{thm_submain}, \ref{thm_main}, and~\ref{thm_third} stating that,
under condition \eqref{costanti},
there exists a minimizer of the Canham--Helfrich functional in 
certain classes 
of generalized Gauss graphs, also enforcing area and enclosed volume constraints, \blu{the latter being the physically relevant setup for biological applications}.
Their proof is a consequence of the direct method in the Calculus of Variations, once lower semicontinuity and compactness are proved.

The main advantage of the generalized Gauss graphs setting is that we are able to cover the physical range \eqref{costanti} for the bending coefficients. One shortcoming is the need of technical conditions to define the classes where the functional is minimized. Nonetheless, regular two-dimensional oriented surfaces always belong to such minimization classes.


The plan of the paper is the following: in Section~\ref{sec:GGG} we present a brief review of generalized Gauss graphs, after which we define the Canham--Helfrich energy of a generalized Gauss graph in Section~\ref{sec:tre}. Section~\ref{sec:lsccpt} is devoted to the main results, Theorems~\ref{thm_submain}, \ref{thm_main}, and~\ref{thm_third}, and is complemented by a regularity result, Theorem~\ref{thm_reg}.

\subsection{Motivations from biological membranes and outlook}
In this section we briefly describe the origin of the Canham--Helfrich functional $\cE$ in \eqref{100} and present an outlook for future research.

In the early 1970's Canham and Helfrich independently proposed a free energy in an effort to model the shape of biological membranes. 
The lipid bilayer that usually constitutes biological membranes is composed of amphiphiles, polar molecules featuring a hydrophilic head and a hydrophobic fatty tail, that are arranged in a fashion so that the tails in the inner part of the bilayer, screened from the watery surrounding environment.

Given the thickness of a few nanometers, one such bilayer can be effectively described as a surface $M$ and the form \eqref{100} of the energy depending \blu{only} on the mean curvature~$H$ of~$M$ responds to the need of explaining the bi-concave shape of red blood cells \cite{Canham1970}.
The competing contribution coming from the Gaussian curvature $K$ was added by Helfrich \cite{Helfrich1973}, whereas the presence of the spontaneous mean curvature $H_0$ takes into account possibly preferred configurations: this is the case in which the asymmetry between the two layers, or the difference in the chemical potential across the membrane determine a natural bending of the membrane, even at rest.

Several derivation for the Canham--Helfrich energy \eqref{100} are available, see \cite{SF2014} and the references therein, which rely on formal expansions of microscopic energies for small thickness.
A more rigorous derivation in terms of $\Gamma$-convergence would be amenable from the variational point of view. Some results in this direction are available. In \cite{PR2009} a complete derivation in dimension two is presented, while the full three-dimensional case is tackled in \cite{LPR2014,LR}, where only partial results are obtained: the $\Gamma-\limsup$ inequality is proved, but the $\Gamma-\liminf$ inequality is proved in the setting of generalized Gauss graphs for a simplified functional. It would be interesting to recover a full $\Gamma$-convergence result also in the three-dimensional case.
\blu{This work sets the stage for possibly tackling this problem in the context of generalized Gauss graphs, especially in light of the sharpness of the bounds on the bending constants.}

%
%
%

\section{Brief theory of generalized Gauss graphs}\label{sec:GGG}
In this section we introduce generalized Gauss graphs and highlight their main properties. We start by introducing some notions from exterior algebra and rectifiable currents.


\subsection{Exterior algebra and rectifiable currents}
We refer the reader to \cite{Federer} for a comprehensive treatise on the theory of currents.

Let $k,N\in\N$ be such that \blu{$1 \leq k \leq N$}.
We define $\bigw^0(\R{N})\coloneqq\R{}$ and we denote by $\bigw^k(\R{N})$ the space of \emph{$k$-covectors} in $\R{N}$, that is the space of $k$-linear alternating forms on $\R{N}$; we denote by $\bigw_k(\R{N})$ the dual space $(\bigw^k(\R{N}))^*=\bigw^k((\R{N})^*)$, called the space of \emph{$k$-vectors} in~$\R{N}$. 
We recall that, if $\{e_1, \ldots, e_N \}$ is a basis of $\R{N}$, then $\{e_{i_1} \wedge \cdots \wedge e_{i_k} : 1 \leq i_1< \cdots<i_k \leq N \}$ is a basis of $\bigw_k(\R{N})$, where $\wedge$ denotes the exterior product.
A $k$-vector $v$ is called a \emph{simple $k$-vector} if it can be written as $v = v_1 \wedge \cdots \wedge v_k$, for some $v_1,\ldots,v_k\in\bigw_1(\R{N})\simeq\R{N}$.

Let $\Omega \subseteq \R{N}$ be an open set.
A \emph{(differential) $k$-form} $\omega$ on $\Omega$ is a map that to each $x \in \Omega$ associates $\omega(x) \in \bigw^k(\R{N})$.
Given $\omega$ a $0$-form on $\Omega$ (that is, a scalar function $\omega\colon\Omega\to\R{}$), we define $\de\omega$ as the $1$-form on $\Omega$ given by the differential of $\omega$; for $k>0$, the definition of the \emph{exterior differential operator} $\de$ is extended from $k$-forms to $(k+1)$-forms through the usual algebra of the exterior product.
We denote by $\cD^k(\Omega)$ the space of $k$-forms with compact support in $\Omega$; the space of \emph{$k$-currents} $\cD_k(\Omega)$ is defined as the dual of $\cD^k(\Omega)$. 
Given a sequence of currents $\{T_n\}_{n \in \N} \subseteq \cD_k(\Omega)$ and a current $T \in \cD_k(\Omega)$, we say that $T_n \wto T$ if and only if $\langle T_n , \omega\rangle \to \langle T, \omega \rangle$ for every $\omega\in\cD^k(\Omega)$, where $\langle \cdot, \cdot \rangle $ denotes the dual product.
We denote by $\partial T \in \mathcal{D}_{k-1}(\Omega)$ the \emph{boundary} of the current $T \in \mathcal{D}_k(\Omega)$, defined as $\langle \partial T, \omega \rangle \coloneqq \langle T, \de \omega \rangle $ for every $\omega \in \mathcal{D}^{k-1}(\Omega)$; we notice that $\de\omega\in\cD^k(\Omega)$ whenever $\omega\in\cD^{k-1}(\Omega)$, that is, exterior differentiation preserves the compactness of the support, so that the duality $\langle\partial T,\omega\rangle$ is well defined.
The \emph{mass} of a current $T \in \mathcal{D}_k(\Omega)$ in the open set $W \subseteq \Omega$ is defined as 
\begin{equation*}
\M{}_W(T)\coloneqq \sup \big\{\langle T, \omega \rangle: \text{$\omega \in \cD^k(W)$, $\|\omega(x) \| \leq 1$ for every $x\in W$}\big\}.
\end{equation*}
Here, $\|\cdot \|$ denotes the comass norm, namely, for $\alpha\in\bigw^k(\R{N})$,
$$\|\alpha\|\coloneqq\sup\{\langle\alpha,v\rangle:\text{$v$ is a simple $k$-vector with $|v|\leq1$}\big\},$$
where $|v|\coloneqq|v_1\wedge\cdots\wedge v_k|$ is the volume of the parallelepiped generated by $v_1,\ldots,v_k$.

Given a set $M \subseteq \R{N}$, we say that $M$ is \emph{$k$-rectifiable} if $M \subseteq \bigcup_{m=0}^{\infty} M_i$, for a certain $\cH^k$-negligible subset $M_0\subseteq\R{N}$ and for certain $k$-dimensional $\cC^1$ surfaces $M_i\subseteq\R{N}$, for $i>0$. 
One can prove that, if $M$ is a $k$-rectifiable set, for $\mathcal{H}^k$-almost every $p \in M$ there exists an approximate tangent space denoted by $T_pM$.
We say that a map $\eta\colon M\to \bigw_k(\R{N})$ is an \emph{orientation} of $M$ if it is $\mathcal{H}^k$-measurable and if $\eta(p)$ is a unit simple $k$-vector that spans the approximate tangent space $T_pM$ for $\mathcal{H}^k$-almost every $p \in M$. We say that a map $\beta\colon M \to \R{}$ is an \emph{integer multiplicity} on $M$ if it is $\mathcal{H}^k$-locally summable and with values in~$\N$.
Finally, $T  \in \mathcal{D}_k(\Omega)$ is a \emph{$k$-rectifiable current with integer multiplicity} if there exist a $k$-rectifiable set  $M\subseteq\R{N}$, an orientation $\eta$ of $M$, and an integer multiplicity $\beta$ on $M$ such that for every $\omega \in \mathcal{D}^k(\Omega)$ we have
\begin{equation*}
\langle T, \omega \rangle =\int_M \langle \omega(p), \eta(p) \rangle \beta(p) \,\de \mathcal{H}^k(p).
\end{equation*}
We denote by $\mathcal{R}_k(\Omega)$ the sets of such currents and write $T= \llb M,\eta,\beta\rrb$. In this case, we have that
\begin{equation}\label{09482}
\M{}_W(T)=\int_{M\cap W} \beta(p)\,\de\cH^{k}(p),
\end{equation}
which simply returns $\cH^{k}(M\cap W)$ if the multiplicity $\beta$ is $1$. 

We now state the celebrated Federer--Fleming compactness theorem, which establishes the compactness result for $k$-rectifiable currents with integer multiplicity.
\begin{theorem}[{\cite[Theorem~4.2.17]{Federer}}]\label{481}
Let $\{T_n\}_{n \in \N}$ be a sequence in $\mathcal{R}_k(\Omega)$ such that $\partial T_n \in \mathcal{R}_{k-1}(\Omega)$ for any $n \in \N$. Assume that for any open set $W$ with compact closure in $\Omega$ there exists a constant $c_{W}>0$ such that
\begin{equation*}
\M{}_W(T_n)+\M{}_W(\partial T_n) < c_W.
\end{equation*}
Then there exist a subsequence $\{n_j\}_{j\in\N}$ and a current $T \in \mathcal{R}_k(\Omega)$ with $\partial T \in \mathcal{R}_{k-1}(\Omega)$ such that  $T_{n_j} \wto T$ as $j \to \infty$.
\end{theorem}


\subsection{Gauss graphs}
We refer the reader to \cite{DoCarmo3d,DoCarmoForme} for the classical notions of differential geometry.

Let $M\subseteq\R{3}$ be a compact two-dimensional manifold of class $\cC^2$; we say that $M$ is \emph{orientable} if there exists a map $\nu\colon M\to\S{2}$ of class $\cC^1$ on $M$ such that, for every $p\in M$, the vector $\nu(p)$ is perpendicular to the tangent space $T_pM$. Once we fix a choice of such a map $\nu$, we say that the manifold $M$ is \emph{oriented} and we call $\nu$ the \emph{Gauss map} of $M$. 
Since $M$ is of class $\cC^2$, the Gauss map is differentiable at any $p \in M$ and, upon identifying $T_{\nu(p)}\S{2}\simeq T_pM$, its differential in $p$, $\de\nu_p\colon T_pM\to T_{\nu(p)}\S{2}$, is a self-adjoint linear operator that has two real eigenvalues $\kappa_1(p)$ and $\kappa_2(p)$, called the \emph{principal curvatures} of~$M$ at $p$. We define the mean and Gaussian curvatures of $M$ at $p$ by
\begin{equation*}
H(p)\coloneqq\kappa_1(p)+\kappa_2(p),\quad K(p)\coloneqq\kappa_1(p)\kappa_2(p).
\end{equation*}
The map $\de\nu_p$ can be extended to a linear map $L_p\colon \R{3}\to\R{3}$ by setting 
\begin{equation}\label{453}
L_p\coloneqq\de\nu_p\circ\bP_p,
\end{equation}
where $\bP_p\colon\R{3}\to T_p M$ denotes the orthogonal projection on the tangent space. Observe that $L_p$ has eigenvalues $\kappa_1(p)$, $\kappa_2(p)$, and $0$; in particular, $\det L_p=0$.

For convenience, we denote by $\R{3}_x$ the space of points $p$ and by $\R{3}_y$ the space where $\nu(p)$ takes its values, and we consider the graph of the Gauss map
\begin{equation}\label{Gg}
G\coloneqq \big\{(p,\nu(p)) \in \R{3}_x \times \R{3}_y : p \in M\big\} \subset \R{3}_x \times \R{3}_y\simeq\R{6}.
\end{equation}
Since $M$ is a two-dimensional manifold of class $\cC^2$, $G$ is a two-dimensional embedded surface in $\R{3}_x \times \R{3}_y$ of class $\cC^1$. We remark that if $M$ has a boundary then also $G$ has a boundary which is given by
$
\partial G=\big\{(p,\nu(p)): p \in\partial M\big\}
$ 
and we notice that if $\partial M = \emptyset$ then $\partial G=\emptyset$.

We now define an orientation on $G$.
We equip $M$ with the orientation induced by $\nu$ and let $\tau(p)\coloneqq*\nu(p)$, where 
\begin{equation*}
*\colon\bigw_1(\R{3}) \to \bigw_2(\R{3})
\end{equation*}
is the Hodge operator.  Notice that $\tau(p) \in \bigw_2(T_pM)$ for every $p \in M$, thus the field $p \mapsto \tau(p)$ is a tangent $2$-vector field on $M$. Then, the action of the Hodge operator is that of, starting from $\nu(p)$, providing the oriented basis vectors of the tangent plane to $M$ at $p$, namely $\tau(p)$.
Let $\Phi\colon M \to M \times \S{2} \subseteq \R{3}_x \times \R{3}_y$ be given by $\Phi(p)\coloneqq (p,\nu(p))$ which is of class $\cC^1$ on $M$. Observe that $G=\Phi(M)$. For each $p \in M$ we have
\begin{align*}
\de \Phi_p\colon T_pM &\to T_pM \times T_{\nu(p)}\S{2}\subseteq \R{3}_x \times \R{3}_y\\
u & \mapsto (u,\de\nu_p(u)).
\end{align*}
Finally, we define $\xi\colon G \to \bigw_2(\R{3}_x \times \R{3}_y)$ as
\begin{equation}\label{def_xi}
\xi(p,\nu(p))\coloneqq  \de\Phi_p(\tau_1(p)) \wedge \de\Phi_p(\tau_2(p)), \quad\text{for}\quad \tau = \tau_1 \wedge \tau_2.
\end{equation}
It is easy to see that $|\xi|\ge 1$, hence we can normalize $\xi$ obtaining 
\begin{equation}\label{483}
\eta\coloneqq\frac{\xi}{|\xi|},
\end{equation}
which is an orientation of $G$.

We now introduce the general setting of generalized Gauss graphs. Let $\{e_1,e_2,e_3\}$ and $\{\varepsilon_1,\varepsilon_2,\varepsilon_3\}$ be the canonical basis of $\R{3}_x$ and $\R{3}_y$, respectively.
Given a $2$-vector $\xi \in \bigw_2(\R{3}_x \times \R{3}_y)$, we define the \emph{stratification} of $\xi$ as the unique decomposition
\begin{equation*}
\xi=\xi_0+\xi_1+\xi_2, \quad \text{where} \quad \xi_0 \in \bigw_2(\R{3}_x), \quad \xi_1 \in\bigw_1(\R{3}_x) \wedge \bigw_1(\R{3}_y), \quad \xi_2 \in \bigw_2(\R{3}_y),
\end{equation*}
given by
\begin{align*}
\xi_0 &=\sum_{1\leq i <j\leq3} \langle \de x_i \wedge \de x_j,\xi\rangle e_i\wedge e_j\eqqcolon\sum_{1\leq i <j\leq3} \xi_0^{ij} e_i\wedge e_j,\\
\xi_1 &=\sum_{1\leq i,j\leq3} \langle \de x_i \wedge \de y_j,\xi\rangle e_i\wedge \varepsilon_j\eqqcolon\sum_{1\leq i,j\leq3} \xi_1^{ij} e_i\wedge \varepsilon_j,\\
\xi_2 &=\sum_{1\leq i <j\leq3} \langle \de y_i \wedge \de y_j,\xi\rangle \varepsilon_i\wedge \varepsilon_j\eqqcolon\sum_{1\leq i <j\leq3} \xi_2^{ij} \varepsilon_i\wedge \varepsilon_j,
\end{align*}
where $\{\de x_1,\de x_2, \de x_3\}$ and $\{\de y_1,\de y_2,\de y_3\}$ denote the dual basis of $\{e_1,e_2,e_3\}$ and $\{\varepsilon_1,\varepsilon_2,\varepsilon_3\}$, respectively. Notice that the three equalities above serve as a definition of $\xi_h^{ij}$; $\xi_0$ and $\xi_2$ are represented by $3\times 3$ skew-symmetric matrices while $\xi_1$ is represented by a $3\times 3$ matrix.

From now on we take $\Omega \subseteq \mathbb R_x^3$ an open set. We indicate by $\curv_2(\Omega)$ the set of the currents $\Sigma=\llb G,\eta,\beta\rrb$ that satisfy
\begin{equation}\label{582}
\begin{cases}
\text{$\Sigma$ and $\partial \Sigma$ are rectifiable currents supported on $\Omega \times \S{2}$,} \\
\displaystyle \langle \Sigma, g\varphi^* \rangle= \int_G g(x,y)|\eta_0(x,y)|\beta(x,y)\, \de \mathcal{H}^2(x,y), \quad\text{for all $g \in \cC_c(\Omega \times \R{3}_y)$,} \\
\langle \partial \Sigma, \varphi \wedge \omega\rangle=0, \quad\text{for all $\omega \in \mathcal{D}^0(\Omega \times \R{3}_y)$,}
\end{cases}
\end{equation}
where 
\begin{equation*}
\varphi(x,y) \coloneqq \sum_{j=1}^3 y_j \de x_j, \quad
\varphi^*(x,y)\coloneqq \sum_{j=1}^3 (-1)^{j+1} y_j \de \hat x_j
\end{equation*}
and $\de \hat x_j=\de x_{j_1} \wedge \de x_{j_2}$ for $1\leq j_1<j_2 \leq 3$, $j_1,j_2 \ne j$.
We can associate with the \blu{regular} Gauss graph $G$ the current $\Sigma_G \in \mathcal{R}_2(\R{3}_x \times \R{3}_y)$ given by $\Sigma_G\coloneqq\llb G,\eta,1\rrb$, and this turns out to be an element of $\curv_2(\Omega)$ (see \cite[Section~2]{AST})\footnote{We notice that the smallest weakly sequentially closed subset of
\begin{equation*}
\{T \in \mathcal{D}_2(\Omega \times \R{3}_y): \M{}_W(T) +\M{}_W(\partial T) <\infty \,\, \forall W \Subset \Omega \times \R{3}\}
\end{equation*}
containing the currents associated with the regular Gauss graphs is a subset of $\curv_2(\Omega)$ (see again \cite[Sec.\,2]{AST}). This shows that $\curv_2(\Omega)$ is a rich enough set.}.
Given $\Sigma=\llb G,\eta,\beta \rrb\in\curv_2(\Omega)$, we let $G^*\coloneqq\{(x,y)\in G: \eta_0(x,y)\neq0\}$ \blu{(notice that $G^*$ is defined only $\cH^2$-a.e.)}.

The geometric meaning of the first condition in \eqref{582} is evident; the one of the second condition is the following: the variable $y$ is orthogonal to the tangent space to $p_1 G$, where we denote by $p_1\colon\R3_x\times\R3_y\to\R3_x$ the projection on the first component; the third condition is the analogue of the second one, on the boundary $\partial\Sigma$.

For a rectifiable current $\Sigma=\llb G,\eta,\beta \rrb$, according with the stratification of $\eta$, we define the strata $\Sigma_i$ by
$$\Sigma_i(\omega)\coloneqq \int_G \langle \omega(x,y),\eta_i(x,y)\rangle \beta(x,y)\, \de \mathcal{H}^2(x,y)$$ 
for every $\omega \in \mathcal{D}^2(\R{3}_x\times\R{3}_y)$.
Given $k\in \{1,2,3 \}$, consider a multi-index $\lambda \in  \{(\lambda_1,\ldots,\lambda_k )\,:\, 0\leq \lambda_1 < \cdots <\lambda_k \leq 2 \}$.
\blu{Letting $|\Sigma_{\lambda_i}|$ be the measure $|\eta_{\lambda_i}|\beta \mathcal{H}^2\res G$,} a generalized Gauss graph $\Sigma \in \curv_2(\Omega)$ is said to be \emph{$\lambda$-special} if 
\begin{equation*}
|\Sigma_{\lambda_i}| \ll |\Sigma_0| \quad \text{for} \quad  i=1,\ldots,k
\end{equation*}
and we write $ \Sigma \in \curv_2^\lambda(\Omega)$. 
We set
$\curv_2^*(\Omega)\coloneqq \curv_2^{(0,1,2)}(\Omega)$
and we call its elements \emph{special generalized Gauss graphs}; in the sequel, we will also make use of the space 
\begin{equation*}
\curv_2^{(0,1)}(\Omega)=\big\{\Sigma\in\curv_2(\Omega): |\Sigma_{1}| \ll |\Sigma_0|\big\}.
\end{equation*}

We introduce the following class of functions.
\begin{defin}\label{stint}
A function $f\colon \Omega \times \R{3}_y\times \big(\bigw_1(\R{3}_x) \wedge \bigw_1(\R{3}_y)\big) \to \R{}$ is said to be a \emph{standard integrand in the setting of $\curv_2(\Omega)$} if
\begin{itemize}
\item[(i)] $f$ is continuous;
\item[(ii)] $f$ is convex in the last variable, \emph{i.e.},
\begin{equation*}
f(x,y,t p+(1-t)q) \leq tf(x,y,p)+(1-t) f(x,y,q),
\end{equation*}
for all $t \in (0,1)$, for all $(x,y) \in \Omega \times \R{3}_y$, and for all $p,q \in \bigw_1(\R{3}_x) \wedge \bigw_1(\R{3}_y)$;
\item[(iii)]$f$ has superlinear growth in the last variable, \emph{i.e.}, there exists a continuous function $\varphi\colon\Omega  \times \R{3}_y \times[0,+\infty) \to [0,+\infty)$, non-decreasing in the last variable and such that $\varphi(x,y,t) \to +\infty$ locally uniformly in $(x,y)$ as $t \to +\infty$, and with
\begin{equation*}
\varphi(x,y,|q|)|q| \leq f(x,y,q)
\end{equation*}
for all $(x,y,q) \in \Omega  \times \R{3}_y \times\big(\bigw_1(\R{3}_x) \wedge \bigw_1(\R{3}_y)\big).$
\end{itemize}
\end{defin}

\begin{remark}
A function $f$ as in Definition \ref{stint} is called (1)-standard integrand in \cite[Definition~3.3]{Del01}.
\end{remark}

\noindent 
The following theorem ensures that an integral functional with a standard integrand as a density is lower semicontinuous.
\begin{theorem}\cite[Theorem~3.2]{Del01}\label{teo_sci}
Let $f$ be a standard integrand in the setting of $\curv_2(\Omega)$ and, for every $\Sigma =\llb G,\eta,\beta \rrb \in \curv_2(\Omega)$, set 
\begin{equation*}
I_f(\Sigma)\coloneqq \int_{G^*} f(x,y,\xi_1(x,y)) |\eta_0(x,y)| \beta(x,y) \,\de \mathcal{H}^2(x,y).
\end{equation*}
Consider a sequence $\{\Sigma_j \}_{j\in\N} \subset \curv_2^{(0,1)}(\Omega)$ such that
\begin{itemize}
\item[(i)] $\Sigma_j \wto \Sigma$, where \blu{$\Sigma \in \mathcal{R}_2(\Omega \times \mathbb{S}^{2})$};
\item[(ii)] $\sup_{j\in\N} I_f(\Sigma_j) <+\infty$.
\end{itemize}
Then
\begin{equation*}
\Sigma \in \curv_2^{(0,1)}(\Omega) \quad \text{and} \quad I_f(\Sigma) \leq \liminf_{j \to \infty} I_f(\Sigma_j).
\end{equation*}
\end{theorem}
\begin{theorem}\cite[Corollary~4.2]{Delladio}\label{teo_comp}
Consider a sequence $\Sigma_j =\llb G_j,\eta_j,\beta_j \rrb\in \curv_2^{*}(\Omega)$ such that
\begin{equation*}
\sup_{j\in\N}\bigg\{ \int_{G_j^*} \!\! \bigg(|(\eta_j)_0(x,y)|+\frac{|(\eta_j)_1(x,y)|^2}{|(\eta_j)_0(x,y)|}+\frac{|(\eta_j)_2(x,y)|^2}{|(\eta_j)_0(x,y)|} \bigg)\beta_j(x,y)\,\de \mathcal{H}^2(x,y)+ \M{}(\partial \Sigma_j)\bigg\} <+\infty.
\end{equation*}
Then there exist a subsequence $\{\Sigma_{j_k}\}_{k\in\N}$ and $\Sigma \in \curv_2^*(\Omega)$ such that $\Sigma_{j_k} \wto \Sigma$ as $k\to\infty$.
\end{theorem}

\section{The Canham--Helfrich energy of a Generalized Gauss Graph}\label{sec:tre}


In this section, we are going to define the Canham--Helfrich energy of a generalized Gauss graph in a way that is the natural extension of the definition for smooth surfaces. Let $H_0\in \R{}$. Here, $M\subseteq\R{3}$ denotes a compact and oriented (with an understood choice of the normal $\nu$) two-dimensional manifold of class~$\cC^2$; recall the definition \eqref{100} of the Canham--Helfrich energy functional $\cE(M)$ on $M$. 

\begin{lemma}[{\cite[Lemma~4.2]{LR}}]
\label{Lemma_4.2_LR}
For $\xi \in \bigw_2(\R{3}_x\times \R{3}_y)$ as in \eqref{def_xi} the following hold true
\begin{equation}\label{3_1}
\begin{cases}
\xi_0=\tau_1 \wedge \tau_2, \\
\xi_1=\tau_1 \wedge \de \nu( \tau_2) -\tau_2 \wedge \de \nu( \tau_1),\\
\xi_1^{ij}=(\tau_1 \otimes \de \nu( \tau_2) -\tau_2 \otimes \de \nu( \tau_1))_{ij},\\
\xi_2=\de\nu (\tau_1) \wedge \de \nu (\tau_2)=\kappa_1 \kappa_2 \tau_1 \wedge \tau_2.\\
\end{cases}
\end{equation}
\end{lemma}

\begin{proof}
By the definition of $\xi$ we have
\begin{equation}\label{formula}
\begin{aligned}
\xi(p,\nu(p))&=(\tau_1(p), \de\nu_p(\tau_1(p))) \wedge (\tau_2(p), \de\nu_p(\tau_2(p)))\\
&=
\tau_1\wedge \tau_2 +\tau_1 \wedge \de \nu_p( \tau_2) -\tau_2 \wedge \de \nu_p( \tau_1)+\de\nu_p (\tau_1) \wedge \de \nu_p (\tau_2).
\end{aligned}
\end{equation}
Then the equalities follow from straightforward computations. 
\end{proof}
\begin{remark}\label{massasotto}
If $M$ is a two-dimensional oriented manifold of class $\cC^2$ with multiplicity $\bar\beta\colon M\to\N$,  $G$ is the Gauss graph associated with $M$ via \eqref{Gg}, and $\Sigma_G\coloneqq \llb G,\eta,\beta\rrb$ with $\beta(x,y)=\bar\beta(x)$, then the equalities
\begin{equation}\label{massaMmult}
\M{}(M)=\int_M \bar\beta(p)\,\de\cH^2(p)=\int_G \frac{\beta(x,y)}{|\xi(x,y)|}\,\de\cH^2(x,y)=\int_G |\eta_0(x,y)|\beta(x,y)\,\de\cH^2(x,y)
\end{equation}
hold true by means of the area formula, \eqref{483}, and the first identity in \eqref{3_1}; here, by $\M{}(M)$ we mean the mass of the current $\llb M,\nu,\bar\beta\rrb$, see \eqref{09482} with $k=2$.
In particular, if $\bar\beta\equiv1$, we obtain
\begin{equation}\label{massaM}
\cH^2(M)=
\int_G |\eta_0(x,y)|\,\de\cH^2(x,y).
\end{equation}

\end{remark}


The next two lemmas are proved in \cite{LR}. We provide the proof in our context for the sake of completeness.
\begin{lemma}[{\cite[Lemma~4.5]{LR}}]
\label{Lemma_4.5_LR}
Let $\Sigma=\llb G,\eta,\beta\rrb \in \curv_2(\Omega)$ be a generalized Gauss graph. Then
\begin{itemize}
\item for $\mathcal{H}^2$-almost every $(x,y) \in G$
\begin{equation}\label{4.27_LR}
\sum_{i=1}^{3} \eta_1^{ij}(x,y) y_i = 0 \quad \text{for all $1 \leq j \leq 3$,}
\end{equation}
\item for $\cH^2$-almost every $(x,y) \in G^*$
\begin{equation}\label{4.28_LR}
\sum_{j=1}^{3} \eta_1^{ij}(x,y) y_j = 0 \quad \text{for all $1 \leq i \leq 3$.}
\end{equation}
\end{itemize}
\end{lemma}

\begin{proof}

As in the proof of \cite[Proposition~2.4]{AST}, we have that
\begin{align*}
\langle \eta(x,y), (y,0) \wedge (0,w) \rangle = 0 \quad \text{for all $w \in \R{3}$ and for $\mathcal{H}^2$-almost every $(x,y) \in G$.}
\end{align*}
From this we deduce that $\sum_{ij} \eta_1^{ij}(x,y) y_i w_j=0$ for all $w \in \R{3}$, which implies \eqref{4.27_LR}.

By \cite[Theorem~2.10(ii)]{AST}, for $\mathcal{H}^2$-almost every $(x,y) \in G^*$, there are an embedded $\cC^1$ surface $S \subset \R{3}$ and a map $\zeta\colon S \to \S{2}$ of class $\cC^1$ such that
\begin{equation*}
\zeta(x)=y, \quad \bigw_2(\Id \oplus \de \zeta_x)(*y)  = \xi(x,y).
\end{equation*}
By Lemma \ref{Lemma_4.2_LR}, 
we obtain, for $i=1,2,3$ and $*y=\tau_1 \wedge \tau_2$,
\begin{equation*}
\sum_{j=1}^3 \xi_1^{ij}y_j=e_i \cdot (\tau_1 \otimes D\zeta(x) \tau_2-\tau_2 \otimes D \zeta(x) \tau_1)y=0,
\end{equation*}
since $D\zeta(x) \tau_k \cdot y=D\zeta (x) \tau_k \cdot \zeta(x)=0$ for $k=1,2$ as $\zeta$ takes values in $\S{2}$. Then \eqref{4.28_LR} is proved recalling \eqref{483}. 
\end{proof}

We recall that the permutation symbols are given by 
$$\eps_{ijk}=\begin{cases}
1 & \text{if $(ijk)$ is an even permutation of $\{1,2,3\}$,} \\
-1 & \text{if $(ijk)$ is an odd permutation of $\{1,2,3\}$,} \\
0 & \text{otherwise.}
\end{cases}$$
For any $z\in\R{3}$, we define
\begin{equation}\label{649}
\Psi_{z}\coloneqq \sum_{i,j,k=1}^3 \varepsilon_{ijk}\,z_k \, \de x_i \wedge \de y_j.
\end{equation}

\begin{lemma}[{\cite[Lemma~4.6]{LR}}]\label{656}
For $L$ as in \eqref{453}, the following formulas hold
\begin{eqnarray}
&&\blu{H=}\tr  L=\nu_1 (\xi_1^{23}-\xi_1^{32}) - \nu_2 (\xi_1^{13}-\xi_1^{31}) + \nu_3 (\xi_1^{12}-\xi_1^{21})= \langle \Psi_{\nu}, \xi_1 \rangle,\label{nuovissimo}\\
&&\blu{K=}\tr  (\cof L)=\nu \cdot (\cof \xi_1) \nu,\nonumber
\end{eqnarray}
where $L$ and $\nu$ are evaluated at $p \in M$ and $\xi$ is evaluated at $(p,\nu(p))$.
\end{lemma}
\begin{proof}

Since $\{\tau_1,\tau_2,\nu\}$ is an orthonormal basis of $\R{3}$, we observe that for any $r \in \R{}$
\begin{align}\label{pol_car}
-r \tr (\cof L)+r^2 \tr L -r^3&=\det (L-r \Id)= \det(\tau_1 | \tau_2 |\nu) \det(L-r\Id)\nonumber\\
&=(L-r\Id)\nu \cdot \left [(L-r\Id)\tau_1 \times (L-r\Id)\tau_2 \right] \\
&=-r(L\tau_1 \times L\tau_2) \cdot \nu + r^2(\tau_1\times L\tau_2-\tau_2 \times L \tau_1)\cdot \nu -r^3 \nonumber,
\end{align}
where we used the fact that $L \nu=0$.
Therefore, from Lemma \ref{Lemma_4.2_LR} we deduce that
\begin{align*}
\tr L &=(\tau_1 \times L\tau_2 -\tau_2 \times L \tau_1) \cdot \nu=\sum_{i,j,k=1}^3 (\tau_{1,i}e_j \cdot L\tau_2-\tau_{2,i}e_j \cdot L \tau_1)\nu_k \varepsilon_{ijk}\\
&=\sum_{i,j,k=1}^3 \xi_1^{ij} \nu_k \varepsilon_{ijk}= \sum_{i<j} \sum_{k=1}^3 (\xi_1^{ij}-\xi_1^{ji}) \nu_k \varepsilon_{ijk} \\
&=\nu_1 (\xi_1^{23}-\xi_1^{32}) - \nu_2 (\xi_1^{13}-\xi_1^{31}) + \nu_3 (\xi_1^{12}-\xi_1^{21}).
\end{align*}
Moreover, from \eqref{3_1} and \eqref{pol_car} we also deduce that
\begin{equation*}
\tr (\cof L)= (L\tau_1 \times L\tau_2) \cdot \nu= (L\tau_1 \wedge L\tau_2) \cdot \xi_0=\kappa_1\kappa_2.
\end{equation*}
Using \eqref{3_1} again and, since $\det(\xi_1)=0$, by \cite[Prop.\,3.21]{Serre}, we have
\begin{equation*}
\nu \cdot \cof(\xi_1)\nu =\nu \cdot \cof(\tau_1 \otimes L\tau_2 -\tau_2 \otimes L\tau_1)\nu = \det(\tau_1 \otimes L\tau_2 -\tau_2 \otimes L\tau_1+ \nu \otimes \nu)\eqqcolon D.
\end{equation*}
We can represent the matrix $\tau_1 \otimes L\tau_2 -\tau_2 \otimes L\tau_1+ \nu \otimes \nu$ with respect to the basis $\{\tau_1, \tau_2,\nu\}$, obtaining 
\begin{align*}
D&= \det \begin{pmatrix} 
L\tau_2 \cdot \tau_1 & L\tau_2 \cdot \tau_2& 0\\
-L\tau_1 \cdot \tau_1 & -L\tau_1 \cdot \tau_2& 0\\
0&0&1
\end{pmatrix} =
\det \begin{pmatrix} 
L\tau_1 \cdot \tau_1 & L\tau_1 \cdot \tau_2\\
L\tau_2 \cdot \tau_1 & L\tau_2 \cdot \tau_2
\end{pmatrix}\\
&=\kappa_1 \kappa_2
\det \begin{pmatrix} 
\tau_1 \cdot \tau_1 & \tau_1 \cdot \tau_2\\
\tau_2 \cdot \tau_1 & \tau_2 \cdot \tau_2
\end{pmatrix}=\kappa_1\kappa_2=\tr \cof L,
\end{align*}
which concludes the proof.
\end{proof}

The next proposition provides the expression of the Canham--Helfrich functional defined on manifolds, seen as regular Gauss graphs. In turns, this suggests how to define the Canham--Helfrich functional for general elements in $\curv_2(\Omega)$.

\begin{proposition}
Fix $y \in \S{2}$ and let
\begin{equation}\label{Xy}
\mathcal{X}_y\coloneqq \bigg \{ \zeta \in \bigw_1(\R{3}_x) \wedge \bigw_1(\R{3}_y) : \sum_{k = 1}^3 \zeta^{k i} y_{k}= \sum_{k=1}^3 \zeta^{i k} y_{k}= \sum_{k=1}^3 \zeta^{k k}= 0 \;\; \text{for all $i=1,2,3$} \bigg \}.
\end{equation}
Let $f_y\colon \mathcal{X}_y \to [0,+\infty)$ be defined by (recall \eqref{649})
\begin{equation}\label{f_y_def}
f_y(\zeta)\coloneqq \alpha_H \langle \Psi_y, \zeta \rangle^2-2\alpha_H H_0 \langle \Psi_y, \zeta \rangle+\alpha_H H_0^2-\alpha_K y \cdot (\cof \zeta) y
\end{equation}
Then, defining $\eta$ as in \eqref{483}, we have
\begin{equation*}
\cE(M)=\int_{\Phi(M)^*}f_y\left (\frac{\eta_1(x,y)}{|\eta_0(x,y)|} \right )|\eta_0(x,y)|\, \de \mathcal{H}^2(x,y).
\end{equation*}
\end{proposition}

\begin{proof}
First observe that, by Lemma~\ref{Lemma_4.5_LR} \blu{and since by \eqref{3_1} the trace of $\xi_1$ is zero}, $\eta_1(x,y)$ belongs to $\mathcal{X}_y$ for almost every $(x,y) \in \Phi(M)^*$.
Moreover, by \eqref{100}, Lemma~\ref{656}, and the area formula, we have
\begin{align*}
\cE(M)&= \int_M \left(\alpha_H (\tr L_p-H_0)^2-\alpha_K\tr(\cof L_p)\right)\de \mathcal{H}^2(p)\\
&=\int_M \Big(\alpha_H\big( \langle \Psi_{\nu(p))}, \xi_1(p,\nu(p)) \rangle-H_0\big)^2-\alpha_K \nu(p) \cdot (\cof \xi_1(p,\nu(p))) \nu(p)\Big)\,\de \mathcal{H}^2(p)\\
&=\int_{\Phi(M)^*} \!\! \bigg(\alpha_H\bigg( \bigg\langle \Psi_{y}, \frac{\eta_1(x,y)}{|\eta_0(x,y)|} \bigg\rangle-H_0\bigg)^2 \! -\alpha_K y \cdot \bigg(\cof \frac{\eta_1(x,y)}{|\eta_0(x,y)|}\bigg) y\bigg) |\eta_0(x,y)|\,\de\mathcal{H}^2(x,y),
\end{align*}
where we have used that $|\xi|=1/|\eta_0|=|\det D\Phi|$.
\end{proof}

We are now ready to define the functional $\cE$ on a generalized Gauss graph.

\begin{defin}
The Canham--Helfrich functional defined on generalized Gauss graphs is the functional $\cE \colon \curv_2(\Omega) \to [-\infty,+\infty]$ defined by 
\begin{equation}\label{698}
\cE(\Sigma)\coloneqq \int_{G^*} f_y\bigg (\frac{\eta_1(x,y)}{|\eta_0(x,y)|} \bigg )|\eta_0(x,y)|\beta(x,y)\,\de\cH^2(x,y),
\end{equation}
for every $\Sigma=\llb G,\eta,\beta\rrb\in\curv_2(\Omega)$.
\end{defin}

%

\section{Existence and regularity of minimizers}\label{sec:lsccpt}

\subsection{Technical lemmas}
For every $\zeta\in\bigw_1(\R{3}_x)\wedge\bigw_1(\R{3}_y)$ and for every $y\in \mathbb S^2$, let us define
\begin{equation}\label{763}
g_y(\zeta)\coloneqq \alpha_H \langle \Psi_y, \zeta \rangle^2-\alpha_K y \cdot (\cof \zeta) y \qquad\text{and}\qquad h_y(\zeta)\coloneqq 2 \alpha_H H_0\langle \Psi_y, \zeta \rangle
\end{equation}
and let us identify $\zeta$ with a vector in $u=u[\zeta]\in\R{9}$ by
\begin{equation*}
u=u[\zeta]\coloneqq (\zeta^{11}, \zeta^{12}, \zeta^{13}, \zeta^{21}, \zeta^{22}, \zeta^{23}, \zeta^{31}, \zeta^{32}, \zeta^{33}).
\end{equation*}
\blu{With these positions, we have (compare with the expression in \eqref{nuovissimo})
$$\langle \Psi_y,\zeta\rangle=(0,y_3,-y_2,-y_3,0,y_1,y_2,-y_1,0)\cdot u=y_1(u_6-u_8)-y_2(u_3-u_7)+y_3(u_2-u_4).$$}
\begin{lemma}\label{770}
Let \eqref{costanti} hold.
The function $g_y\colon \bigw_1(\R{3}_x)\wedge\bigw_1(\R{3}_y)\to\R{}$ defined in \eqref{763} is represented by a quadratic form $u \mapsto u \cdot A_y u$ on $\R{9}$, where
\footnotesize{
\begin{equation*} 
A_y= 
\begin{pmatrix} 
0 & 0 & 0 & 0 & -\frac{\alpha_K}{2}y_3^2 & \frac{\alpha_K}{2} y_2 y_3 & 0 & \frac{\alpha_K}{2}y_2y_3  & -\frac{\alpha_K}{2}y_2^2 \\
0 & \alpha_Hy_3^2 & -\alpha_Hy_2y_3 &-\gamma y_3^2 & 0 &\gamma y_1y_3& \gamma y_2y_3 & -\alpha_Hy_1y_3 & \frac{\alpha_K}{2}y_1y_2\\
0 & -\alpha_Hy_2y_3 & \alpha_Hy_2^2 &\gamma y_2y_3 & \frac{\alpha_K}{2} y_1 y_3 & \blu{-}\alpha_H y_1 y_2 & -\gamma y_2^2 &\gamma y_1 y_2  & 0\\
0 &-\gamma y_3^2 &\gamma y_2 y_3 &  \alpha_H y_3^2 & 0 & -\alpha_H y_1 y_3 & -\alpha_H y_2 y_3 &\gamma y_1y_3 & \frac{\alpha_K}{2} y_1 y_2 \\
-\frac{\alpha_K}{2} y_3^2 & 0 & \frac{\alpha_K}{2} y_1 y_3 & 0 & 0 & 0 & \frac{\alpha_K}{2}y_1y_3 & 0 & -\frac{\alpha_K}{2} y_1^2\\
\frac{\alpha_K}{2} y_2 y_3 &\gamma y_1 y_3  & -\alpha_H y_1 y_2 & -\alpha_Hy_1 y_3 & 0 & \alpha_H y_1^2 &\gamma y_1y_2 &-\gamma y_1^2 & 0\\
0 &\gamma y_2y_3 &-\gamma y_2^2 & -\alpha_H y_2 y_3 & \frac{\alpha_K}{2} y_1 y_3 &\gamma y_1 y_2 & \alpha_H y_2^2 & -\alpha_H y_1 y_2 & 0\\
\frac{\alpha_K}{2} y_2 y_3 & -\alpha_H y_1 y_3 &\gamma y_1 y_2 &\gamma y_1y_3 & 0 &-\gamma y_1^2 & -\alpha_H y_1y_2 & \alpha_H y_1^2 & 0\\
-\frac{\alpha_K}{2} y_2^2 & \frac{\alpha_K}{2} y_1 y_2 & 0 & \frac{\alpha_K}{2} y_1y_2 & -\frac{\alpha_K}{2} y_1^2 & 0 & 0 & 0 & 0
\end{pmatrix}
\end{equation*}}
\normalsize
for $\gamma \coloneqq \alpha_H - \alpha_K/2$. Let
\begin{equation*}
\begin{split}
v(-\alpha_K/2)\coloneqq &\;
\begin{pmatrix}
y_1^2-1\\
y_1 y_2\\
y_1 y_3\\
y_1 y_2\\
y_2^2-1\\
y_2 y_3\\
y_1 y_3\\
y_2 y_3\\
y_3^2 -1
\end{pmatrix},\qquad
v(2\alpha_H -\alpha_K /2)\coloneqq
\begin{pmatrix}
     0\\
-y_3\\
 y_2\\
 y_3\\
     0\\
    -y_1\\
-y_2\\
     y_1\\
     0
\end{pmatrix}, \\
v_1(\alpha_K/2)\coloneqq &\;
\begin{pmatrix}
 2 y_1y_2y_3\\
 y_3 y_2^2 -y_3 y_1^2\\ 
y_2y_3^2-y_2\\ 
 y_3 y_2^2 -y_3 y_1^2\\    
 -2y_1 y_2y_3\\                     
 y_1-y_1 y_3^2\\                 
y_2y_3^2-y_2\\                     
 y_1-y_1 y_3^2\\ 
0
\end{pmatrix}, \qquad
v_2(\alpha_K/2)\coloneqq
\begin{pmatrix}
y_1y_2^2-y_1 y_3^2 \\ 
 y_2^3-y_2  \\                 
y_3y_1^2+ y_3y_2^2\\ 
 y_2^3-y_2\\ 
y_1-y_1y_2^2\\                      
0\\                 
y_3y_1^2+ y_3y_2^2\\                      
0\\ 
-y_1^3-y_1y_2^2
\end{pmatrix}. 
\end{split}
\end{equation*}
\blu{Then these vectors are eigenvectors of the matrix $A_y$ with corresponding eigenvalues $-\alpha_K /2 $, $2\alpha_H -\alpha_K /2$, and $\alpha_K /2$ with multiplicities $1$, $1$, and $2$, respectively.
The six vectors
\begin{equation*}
\begin{split}
v_1(0)\coloneqq &\,
\begin{pmatrix}
y\\
0\\
0
\end{pmatrix},\;\;
v_2(0)\coloneqq
\begin{pmatrix}
0\\
y\\
0
\end{pmatrix},\;\;
v_3(0)\coloneqq
\begin{pmatrix}
0\\
0\\
y
\end{pmatrix}, \\
v_4(0)\coloneqq &\,
\begin{pmatrix}
y_1e_1\\
y_2e_1\\
y_3e_1
\end{pmatrix}, \;\;
v_5(0)\coloneqq 
\begin{pmatrix}
y_1e_2\\
y_2e_2\\
y_3e_2
\end{pmatrix}, \;\;
 v_6(0)\coloneqq 
 \begin{pmatrix}
y_1e_3\\
y_2e_3\\
y_3e_3
\end{pmatrix}
\end{split}
\end{equation*}
generate the $5$-dimensional subspace associated with the eigenvector $0$.}

The function $h_y\colon \bigw_1(\R{3}_x)\wedge\bigw_1(\R{3}_y)\to\R{}$ defined in \eqref{763}  is represented by a linear map $u \mapsto u \cdot v_y$ where
$ 
v_y\coloneqq \blu{-}2\alpha_H H_0 v(2\alpha_H-\alpha_K/2). 
$ 

Moreover, we have that
\begin{equation}\label{1000}
\begin{split}
&\, \spann \{v_1(0),v_2(0),v_3(0),v_4(0),v_5(0),\blu{v_6(0)},v(-\alpha_K/2)\} \\
=&\, \spann \left\{
\begin{pmatrix}y\\0\\0\end{pmatrix},
\begin{pmatrix}0\\y\\0\end{pmatrix},
\begin{pmatrix}0\\0\\y\end{pmatrix},
\begin{pmatrix}y_1e_1\\y_2e_1\\y_3e_1\end{pmatrix},
\begin{pmatrix}y_1e_2\\y_2e_2\\y_3e_2\end{pmatrix},
\begin{pmatrix}y_1e_3\\y_2e_3\\y_3e_3\end{pmatrix},
\begin{pmatrix}e_1\\e_2\\e_3\end{pmatrix}
\right\} 
\end{split}
\end{equation}
\normalsize
and by the isomorphism $\zeta \mapsto u[\zeta]$ the space $\mathcal{X}_y$ introduced in~\eqref{Xy} transforms to
\begin{equation}\label{tildeXy}
\widetilde\cX_y\coloneqq \left \{ u \in \R{9} : u \perp  \spann \{v_1(0),v_2(0),v_3(0),v_4(0),v_5(0),\blu{v_6(0)},v(-\alpha_K/2)\} \right \}.
\end{equation}
\end{lemma}
\begin{proof}
The claims follow by straightforward calculations. 
\blu{To prove} \eqref{1000}, we observe that,
\begin{align*}
v(-\alpha_K/2)=&\, y_1 \begin{pmatrix}y\\0\\0\end{pmatrix}+y_2\begin{pmatrix}0\\y\\0\end{pmatrix}+y_3\begin{pmatrix}0\\0\\y\end{pmatrix}-\begin{pmatrix}e_1\\e_2\\e_3\end{pmatrix}
\end{align*}
and this concludes the proof.
\end{proof}

Lemma~\ref{770} shows that the quadratic form~$A_y$ (and therefore the function~$g_y$) has both a negative eigenvalue and the zero eigenvalue, which prevent positive definiteness. 
Nonetheless, since the space $\cX_y$ defined in \eqref{Xy} transforms to $\widetilde\cX_y$ defined in \eqref{tildeXy}, which is the orthogonal to the directions where there is loss of positive definiteness, we are able to prove, in the next Proposition, that it is possible to modify the integrand~$f_y$ defined in~\eqref{f_y_def} to obtain the new function~$\tilde f$ defined in~\eqref{effepiccolotilde} below, which  is a standard integrand in the sense of Definition~\ref{stint}.

\begin{proposition}\label{1035}
Let \eqref{costanti} hold. For $y \in \S{2}$, define the map $F_y\colon\R{9}\to \R{}$ 
\begin{align*}
 F_y(u)&\coloneqq g_y(u)-h_y(u)+\alpha_H H_0^2+ \frac{\alpha_K}{2}|\pi_0 u|^2 +\alpha_K |\pi_{-\alpha_K/2}u|^2\\
&=u \cdot  A_y u  -u \cdot v_y+\alpha_H H_0^2+ \frac{\alpha_K}{2}|\pi_0 u|^2 +\alpha_K |\pi_{-\alpha_K/2}u|^2,
\end{align*}
where $g_y,h_y$ are defined as in \eqref{763}, $\pi_0,\pi_{-\alpha_K/2}\colon\R9\to\R9$ are the orthogonal projections on $\spann\{v_1(0),\ldots,v_5(0)\}$ and $\spann\{v(-\alpha_K/2)\}$, respectively.
Moreover, let
\begin{equation}\label{effepiccolotilde}
\tilde f\colon \Omega \times \S{2} \times \left ( \bigw_1(\R{3}_x) \wedge \bigw_1(\R{3}_y)\right )\to \R{}, 
\qquad \tilde f(x,y,\zeta)\coloneqq F_y(u[\zeta]).
\end{equation}
 Then $\tilde f$ is continuous, convex in the third variable, and there exist two constants $c_1>0$ and $c_2\geq0$ 
such that
\begin{equation}\label{1045}
\tilde f(x,y,\zeta) \geq 
c_1|\zeta|^2-c_2.
\end{equation}
In particular, $\tilde f$ has uniform superlinear growth in the third variable.
\end{proposition}

\begin{proof}
Let $\pi_{2\alpha_H-\alpha_K/2},\pi_{\alpha_K/2}\colon\R9\to\R9$ be the orthogonal projections on $\spann\{v(2\alpha_H-\alpha_K/2)\}$ and $\spann\{v_1(\alpha_K/2),v_2(\alpha_K/2)\}$, respectively.
For every $u \in \R{9}$, by Lemma~\ref{770}, we have
\begin{equation}\label{quelleuguaglianzesotto}
\begin{split}
F_y(u)=&\;-\frac{\alpha_K}{2}|\pi_{-\alpha_K/2}u|^2+\frac{\alpha_K}{2}|\pi_{\alpha_K/2}u|^2+\Big(2\alpha_H-\frac{\alpha_K}{2}\Big) |\pi_{2\alpha_H-\alpha_K/2}u|^2\\
&\;\blu{+}2\alpha_H H_0\, u\cdot v(2\alpha_H-\alpha_K/2)+ \frac{\alpha_K}{2}|\pi_0 u|^2 +\alpha_K |\pi_{-\alpha_K/2}u|^2+\alpha_H H_0^2\\
=&\;\frac{\alpha_K}{2}|\pi_{-\alpha_K/2}u|^2+\frac{\alpha_K}{2}|\pi_{\alpha_K/2}u|^2+
\Big(2\alpha_H-\frac{\alpha_K}{2}\Big) |\pi_{2\alpha_H-\alpha_K/2}u|^2\\
&\;+ \frac{\alpha_K}{2}|\pi_0 u|^2\blu{+}2\alpha_H H_0\, u\cdot v(2\alpha_H-\alpha_K/2)+\alpha_H H_0^2.
\end{split}
\end{equation}
By \eqref{costanti}, we deduce that $F_y$ is convex (and therefore continuous) in $u$, so that $\tilde f$ is convex (and therefore continuous) in the third variable. 
Moreover, by reconstructing the norm $|u[\zeta]|^2=|\zeta|^2$ from the projections $\pi_\bullet$ and by recalling that they are $1$-Lipschitz functions, we have that
\begin{equation*}
\tilde f(x,y,\zeta)=F_y(u[\zeta]) \geq \min \left \{\frac{\alpha_K}{2},2\alpha_H-\frac{\alpha_K}{2}\right \}|\zeta|^2-2\sqrt2\alpha_H |H_0| |\zeta|+\alpha_H H_0^2
\end{equation*}
(the factor $\sqrt2=|v(2\alpha_H-\alpha_K/2)|$ comes from Schwarz inequality),
from which we deduce the boundedness from below of $\tilde f$ and \eqref{1045}, with (a possible choice of)
\begin{equation*}
c_1=\frac{1}{4}\min\{\alpha_K,4\alpha_H-\alpha_K\} \quad \text{and} \quad c_2=\alpha_HH_0^2\bigg(\frac{8\alpha_H}{\min\{\alpha_K,4\alpha_H-\alpha_K\}}-1\bigg).
\end{equation*}
Finally, the continuity of $\tilde f$ with respect to $y$ follows from the structure of the matrix $A_y$ and of the vector $v_y$ in Lemma~\ref{770}.
\end{proof}

\begin{proposition}\label{tildemeglio}
Let $\tilde f$ be the function defined in \eqref{effepiccolotilde}. 
The, for every $\Sigma=\llb G,\eta,\beta\rrb\in\curv_2(\Omega)$, it holds that
\begin{equation}\label{1108}
\cE(\Sigma)= \int_{G^*}\tilde f\bigg (x,y,\frac{\eta_1(x,y)}{|\eta_0(x,y)|} \bigg )|\eta_0(x,y)|\beta(x,y)\,\de\cH^2(x,y).
\end{equation}
%
%
\end{proposition}
\begin{proof}
Let $(x,y)\in G^*$. 
By \blu{Lemma~\ref{Lemma_4.5_LR}} and Lemma~\ref{770} we have that $u[\xi_1(x,y)]\in\widetilde\cX_y$, from which we obtain that
$\pi_0u[\xi_1(x,y)]=
\pi_{-\alpha_K/2}u[\xi_1(x,y)]=0$.
Keeping \eqref{f_y_def}, \eqref{763}, and \eqref{effepiccolotilde} into account, this implies that 
$$\tilde f(x,y,\xi_1(x,y))=F_y(u[\xi_1(x,y)])=f_y(\xi_1(x,y)),$$
which, by \eqref{698}, implies \eqref{1108}. 
\end{proof}

\begin{lemma}\label{1263}
Let $A\Subset\Omega$ and let $\Sigma_j=\llb G_j,\eta_j,\beta_j\rrb\in\curv_2(\Omega)$ be such that $\spt\Sigma_j\subseteq A\times\S2$ for every $j\in\N$ and $\Sigma_j\wto\Sigma=\llb G,\eta,\beta\rrb\in\curv_2(\Omega)$ as $j\to\infty$. Then $\spt\Sigma\subseteq A\times\S2$ and 
\begin{equation}\label{1265}
\lim_{j\to\infty}\int_{G_j} |(\eta_j)_0(x,y)|\beta_j(x,y)\,\de\cH^{2}(x,y)=\int_G |\eta_0(x,y)|\beta(x,y)\,\de\cH^{2}(x,y).
\end{equation}
In particular, if $M_j,M$ are two-dimensional oriented manifold of class $\cC^2$ contained in $A$, if $G_j,G$ are the associated Gauss graphs by \eqref{Gg}, and $\Sigma_j=\Sigma_{G_j}=\llb G_j,\eta_j,1\rrb$, $\Sigma=\Sigma_G=\llb G,\eta,1\rrb$ are the associated currents, if $\Sigma_j\wto\Sigma$, then $\cH^2(M_j)\to\cH^2(M)$.
\end{lemma}
\begin{proof}
We first observe that the condition on the supports is closed, so that $\spt\Sigma\subseteq A\times\S2$.
Let $g\in \cC_c(\Omega\times\R3_y)$ be such that $g=1$ on $A\times\S2$. Then the convergence
$$\int_{G_j} |(\eta_j)_0(x,y)|\beta_j(x,y)\,\de\cH^2(x,y)=\Sigma_j(g\varphi^*)\to\Sigma(g\varphi^*)=\int_G |\eta_0(x,y)|\beta(x,y)\,\de\cH^{2}(x,y)$$
follows immediately by \eqref{582}.
The proof of the last statement 
is obtained by combining \eqref{massaM} and \eqref{1265}:
$$\lim_{j\to\infty} \cH^2(M_j)=\lim_{j\to\infty}\int_{G_j} |(\eta_j)_0(x,y)|\,\de\cH^{2}(x,y)=\int_G |\eta_0(x,y)|\,\de\cH^{2}(x,y)=\cH^2(M).$$
This concludes the proof.
\end{proof}


\subsection{Minimization problems}
In this section we study various minimization problems for the energy $\cE$ in \eqref{698}.
In the first two (see Theorems~\ref{thm_submain} and~\ref{thm_main} below), reasonable sufficient conditions for unconstrained minimization are provided. In the third one (see Theorem~\ref{thm_third} below), we tackle constrained minimization in terms of prescribed enclosed volume and surface area for a closed membrane.

For $A\Subset\Omega$ and $c>0$, 
we define the class
\begin{equation}\label{1264}
\cX_{A,c}^{(0,1)}(\Omega)\coloneqq \big\{ \Sigma=\llb G,\eta,\beta\rrb\in\curv_2^{(0,1)}(\Omega) : \spt\Sigma\subseteq A\times\S2,\, \M{}(\partial\Sigma)+\M{}(\Sigma)\leq c\big\}
\end{equation}
of generalized Gauss graphs with compact support and equi-bounded masses.
Our first existence result is the following.
\begin{theorem}\label{thm_submain}
Let \eqref{costanti} hold. The minimization problem
\begin{equation}\label{1286}
\min\Big\{\cE(\Sigma) : \Sigma\in \cX_{A,c}^{(0,1)}(\Omega)\Big\}
\end{equation}
has a solution.
\end{theorem}
\begin{proof}
Let $c_2$ be the constant in \eqref{1045} and, for every $\Sigma=\llb G,\eta,\beta\rrb\in\curv_2^{(0,1)}(\Omega)$, define the functional 
\[\begin{split}
\cE^{(0,1)}(\Sigma)\coloneqq& \int_{G^*} \left(\tilde f \left (x,y,\frac{\eta_1(x,y)}{|\eta_0(x,y)|} \right)+c_2 \right)|\eta_0(x,y)| \beta(x,y)\, \de \mathcal{H}^2(x,y) \\
=&\, \cE(\Sigma)+c_2\int_{G^*} |\eta_0(x,y)| \beta(x,y)\, \de\cH^2(x,y),
\end{split}\]
where the last equality follows from Proposition~\ref{tildemeglio}. 
Inequality \eqref{1045} allows us to apply Theorem~\ref{teo_sci} and obtain that $\cE^{(0,1)}$ is lower semicontinuous in $\curv_2^{(0,1)}(\Omega)$.
By Lemma~\ref{1263}, it follows that also the functional $\cE$ is lower semicontinuous in $\curv_2^{(0,1)}(\Omega)$.
By Theorems~\ref{481} and~\ref{teo_sci}, 
any minimizing sequence $\Sigma_j=\llb G_j,\eta_j,\beta_j\rrb\in\cX_{A,c}^{(0,1)}(\Omega)$ for $\cE$ is compact in $\cX_{A,c}^{(0,1)}(\Omega)$. The thesis then follows from the direct method of the Calculus of Variations.
\end{proof}

\blu{Inequality \eqref{1045} and Lemma~\ref{1263} suggest that it is not necessary to bound the entire $\int_{G^*} \frac{\beta}{|\eta_0|} \de\mathcal{H}^2$ for $\Sigma=\llb G,\eta,\beta\rrb \in \curv_2^*(\Omega)$ in order to apply Theorem~\ref{teo_comp}, so that we can consider the class}

\begin{equation}\label{1322}
\begin{split}
\cX_{A,c}^{*}(\Omega)\coloneqq\bigg\{&\,\Sigma=\llb G,\eta,\beta\rrb\in\curv_2^{*}(\Omega) : \spt\Sigma\subseteq A\times\S2,\\ 
&\, \M{}(\partial\Sigma)+ \int_{ G^*}  \bigg(|\eta_0(x,y)|+\frac{|\eta_2(x,y)|^2}{|\eta_0(x,y)|}\bigg)\beta(x,y)\,\de\cH^2(x,y)\leq c \bigg\}.
\end{split}
\end{equation}
\blu{The bound on  $\int_{G^*}\frac{|\eta_1(x,y)|^2}{|\eta_0(x,y)|^2}|\eta_0(x,y)| \beta(x,y)\, \de \mathcal{H}^2(x,y)$, together with the one on the second term in \eqref{1322}, imply the boundedness of the mass of $\Sigma$. Moreover, these bounds are needed in order to have closedness in the class $\curv_2^*(\Omega)$, which in general is not closed, contrary to $\curv_2(\Omega)$. In particular, for the regular Gauss graph $G$ of a manifold $M$, they imply an $L^4$-bound on the curvatures of $M$, since
\begin{equation*}
\int_{G^*} \left (|\eta_0|+\frac{|\eta_1|^2}{|\eta_0|}+\frac{|\eta_2|^2}{|\eta_0|}\right ) \de \mathcal{H}^2 = \int_M |\xi|^2 \de \mathcal{H}^2 = \int_M \left (H(x)^2+ (1-K(x))^2\right )\, \de \mathcal{H}^2(x),
\end{equation*}
for the proof see \cite[Proposition~1.1 and Example~1.2]{AST}.}
We present now our second existence result.
\begin{theorem}\label{thm_main}
Let \eqref{costanti} hold. The minimization problem
\begin{equation}\label{1343}
\min\Big\{\cE(\Sigma) : \Sigma\in \cX_{A,c}^{*}(\Omega)\Big\}
\end{equation}
has a solution.
\end{theorem}
\begin{proof}
Let us consider a minimizing sequence $\Sigma_j=\llb G_j,\eta_j,\beta_j\rrb\in\cX_{A,c}^{*}(\Omega)$ for the functional~$\cE$. 
By Proposition~\ref{tildemeglio} and \eqref{1045}, we obtain
\[\begin{split}
\cE(\Sigma_j)=& \int_{G_j^*}\tilde f \left (x,y,\frac{(\eta_j)_1(x,y)}{|(\eta_j)_0(x,y)|} \right )|(\eta_j)_0(x,y)|\beta_j(x,y)\, \de \mathcal{H}^2(x,y) \\
\geq&\, c_1\int_{G_j^*}\frac{|(\eta_j)_1(x,y)|^2}{|(\eta_j)_0(x,y)|^2}|(\eta_j)_0(x,y)| \beta_j(x,y)\, \de \mathcal{H}^2(x,y) \\
&\, -c_2\int_{G_j^*}|(\eta_j)_0(x,y)| \beta_j(x,y)\, \de \mathcal{H}^2(x,y).
\end{split}\]
Now, by \eqref{1322}, the minimizing sequence satisfies the hypotheses of Theorem \ref{teo_comp} and therefore there exist a subsequence $\{\Sigma_{j_k}\}_{k\in\N}$ and a special generalized Gauss graph $\Sigma_\infty \in \curv_2^{*}(\Omega)$ such that $\Sigma_{j_k} \wto \Sigma_\infty$ as $k\to\infty$. The thesis follows from the direct method of the Calculus of Variations.
\end{proof}

\begin{remark}
We called the minimization problems \eqref{1286} and \eqref{1343} \emph{unconstrained} because the classes $\cX_{A,c}^{(0,1)}(\Omega)$ in \eqref{1264} and $\cX_{A,c}^*(\Omega)$ in \eqref{1322} do not contain geometric constraints, namely, there are no generalized Gauss graphs excluded from these classes based on their geometry.
In particular, this allows us to consider the zero current $\Sigma=0$ as a competitor for both minimization problems, and it turns out to be an absolute minimizer if $H_0=0$. 
Indeed, in this case, \eqref{1045} becomes $\tilde f(x,y,\zeta)\geq c_1|\zeta|^2$, so that $\cE\geq0$. 
Notice that also a generalized Gauss graph $\Pi$ supported on a plane ($H=K=0$) has zero energy, showing that both \eqref{1286} and \eqref{1343} have no unique solution.

On the other hand, if $H_0\neq0$, observe that a sphere $\Sigma$ (or a portion of it, compatibly with~$A$) with mean curvature $H=H_0$ makes the functional~$\cE$ negative. Indeed, since for spheres there holds $K=H^2/4$, we have $\cE(\Sigma)=-\alpha_K H_0^2 \cH^2(\Sigma)/4<0=\cE(0)<\alpha_HH_0^2=\cE(\Pi)$.
\end{remark}

Given $\Sigma=\llb G,\eta,\beta\rrb\in\curv_2(\Omega)$, we define
\begin{equation}\label{1442}
\cA(\Sigma)\coloneqq  \int_G |\eta_0(x,y)|\beta(x,y)\,\de\cH^{2}(x,y).
\end{equation}
In light of Remark~\ref{massasotto}, if $\Sigma$ is a regular Gauss graph with multiplicity, the quantity $\cA(\Sigma)$ has the geometric interpretation of mass of $p_1\Sigma$, see \eqref{massaMmult}; in particular, if $\beta\equiv1$, then $\cA(\Sigma)=\cH^2(M)$, the \emph{area} of the manifold $M\coloneqq p_1\Sigma$, see \eqref{massaM}.

We also define the quantity
\begin{equation}\label{1451}
\cV(\Sigma)\coloneqq\frac13\int_G (x\cdot y)\, |\eta_0(x,y)|\beta(x,y)\,\de\cH^2(x,y).
\end{equation}
If $\Sigma$ is a closed ($\partial\Sigma=0$) regular Gauss graph with multiplicity $\beta\equiv1$, by a simple application of the Divergence Theorem, the quantity $\cV(\Sigma)$ has the geometric interpretation of the \emph{enclosed volume} in $M\coloneqq p_1\Sigma$. 
\blu{Indeed, if $M=\partial A$, then by means of the area formula we get
\begin{equation*}
\frac13\int_G (x\cdot y)\, |\eta_0(x,y)|\,\de\cH^2(x,y)=\frac13\int_M p\cdot \nu(p)\,\de\cH^2(p)=\frac13\int_A \div(p)\, \de p =\cH^3(A).
\end{equation*}}

\begin{lemma}\label{1263v}
Let $A\Subset\Omega$ and let $\Sigma_j=\llb G_j,\eta_j,\beta_j\rrb\in\curv_2(\Omega)$ be such that $\spt\Sigma_j\subseteq A\times\S2$ and $\partial\Sigma_j=0$ for every $j\in\N$ and $\Sigma_j\wto\Sigma=\llb G,\eta,\beta\rrb\in\curv_2(\Omega)$ as $j\to\infty$. Then $\spt\Sigma\subseteq A\times\S2$, $\partial\Sigma=0$, and 
\begin{equation}\label{12652}
\lim_{j\to\infty}\int_{G_j} (x\cdot y)\, |(\eta_j)_0(x,y)|\beta_j(x,y)\,\de\cH^{2}(x,y)=\int_G (x\cdot y)\, |\eta_0(x,y)|\beta(x,y)\,\de\cH^{2}(x,y).
\end{equation}
In particular, if $M_j=\partial E_j$ and $M=\partial E$ for $E_j,E$ sets of class $\cC^2$ contained in $A$, if $G_j,G$ are the associated Gauss graphs by \eqref{Gg}, and $\Sigma_j=\Sigma_{G_j}=\llb G_j,\eta_j,1\rrb$, $\Sigma=\Sigma_G=\llb G,\eta,1\rrb$ are the associated currents, if $\Sigma_j\wto\Sigma$, then $\cH^3(E_j)\to\cH^3(E)$.
\end{lemma}
\begin{proof}
The proof is the same as that of Lemma~\ref{1263}.
\end{proof}

Next we study constrained minimization problems, namely we prescribe the surface area and the enclosed volume. Given $a,v>0$, we define the class\blu{es
\begin{equation*}\label{nuovaclasse}
\begin{split}
\cX_{A,c;a,v}^{(0,1)}(\Omega)\coloneqq\bigg\{&\,\Sigma=\llb G,\eta,\beta\rrb\in\curv_2^{(0,1)}(\Omega) : \spt\Sigma\subseteq A\times\S2, \partial\Sigma=0,\\ 
&\, \M{}(\Sigma) \leq c, \cA(\Sigma)=a,\cV(\Sigma)=v \bigg\}.
\end{split}
\end{equation*}}
\begin{equation*}\label{1467}
\begin{split}
\cX_{A,c;a,v}^{*}(\Omega)\coloneqq\bigg\{&\,\Sigma=\llb G,\eta,\beta\rrb\in\curv_2^{*}(\Omega) : \spt\Sigma\subseteq A\times\S2, \partial\Sigma=0,\\ 
&\, \int_{ G^*}  \frac{|\eta_2(x,y)|^2}{|\eta_0(x,y)|}\beta(x,y)\,\de\cH^2(x,y)\leq c, \cA(\Sigma)=a,\cV(\Sigma)=v \bigg\}.
\end{split}
\end{equation*}
\blu{In order for two-dimensional closed oriented manifolds of class $\cC^2$ to belong to these classes, we enforce the isoperimetric inequality}
\begin{equation}\label{iso}
36\pi\,v^2\leq a^3.
\end{equation}
\begin{theorem}\label{thm_third}
Let \eqref{costanti} hold and let $a,v>0$ satisfy \eqref{iso}.
The minimization problems
\begin{equation}\label{1480}
\blu{\min\Big\{\cE(\Sigma) : \Sigma\in \cX_{A,c;a,v}^{(0,1)}(\Omega)\Big\}
,} \quad \min\Big\{\cE(\Sigma) : \Sigma\in \cX_{A,c;a,v}^{*}(\Omega)\Big\}
\end{equation}
have a solution.
\end{theorem}
\begin{proof}
The proof is the same as that of Theorems~\ref{thm_submain} and \ref{thm_main}, upon noting that Lemmas~\ref{1263} and~\ref{1263v} provide the continuity for the area and enclosed volume constraints.
\end{proof}

We conclude this subsection with two remarks on the necessity of assumption \eqref{costanti}.

\begin{remark}[$4\alpha_H \leq \alpha_K$]\label{tildemegliodavvero}
In the case, then there exists a constant $r\geq 0$ such that $\alpha_K=4\alpha_H+r$. 
For the Gauss graph $G$ of a smooth surface $M$, we have
\begin{equation*}
    \mathcal{H}^2(G)= \int_M |\xi(p,\nu(p))| \,\de \mathcal{H}^{2}(p)=\int_M \sqrt{4H(p)^2+ (1-K(p))^2}\, \de \mathcal{H}^2(p),
\end{equation*}
where $\xi$ is defined in \eqref{def_xi}.
We consider $ M_j=\partial B_{1/j}$, where $B_{1/j}$ is the ball of radius $1/j$ centered in the origin, and we let $\Sigma_{\blu{j}}\coloneqq \Sigma_{G_j}=\llb G_j,\eta_j,1\rrb$.
Since the principal curvatures of $M_j$ are both equal to $j$, we get from the above formula
\begin{align*}
    \mathbb{M}(\Sigma_j)=\mathcal{H}^2(G_j) \leq \frac{4\pi}{j^2}\sqrt{j^4+14j^2+1},
\end{align*}
which is uniformly bounded for every $j \in \N\setminus\{0\}.$
Thus, for $\Omega=B_2$, we have that $\Sigma_j\in \curv_2^{\blu{(0,1)}}(\Omega)$ for every $j\in\N$ and, since $\partial \Sigma_j=0$, we also have that $\Sigma_j$ belongs  to $\cX_{A,c}^{(0,1)}(\Omega)$ 
for every $j\in\N\setminus\{0\}$,  for a suitable choice of $A$ and $c$.
Since $\Sigma_j$ is a regular Gauss graph, $\cE(\Sigma_j)=\cE(M_j)$, so that, using the expression in \eqref{100}, we obtain
\begin{equation}\label{1434}
\cE(M_j)=4\pi \bigg(\frac{\alpha_HH_0^2}{j^2}-\frac{4\alpha_HH_0}{j}-r\bigg),
\end{equation}
using the fact that $H^2=4K$ for spheres.
We now consider two cases.
\begin{enumerate}
\item[(1)] $r>0$: the functional $\cE$ is no longer lower semicontinuous, since $\Sigma_j \wto 0$ and, by \eqref{1434}, $\displaystyle\liminf_{j\to\infty}\cE(M_j)=-4\pi r < 0= \cE(0)$. 
\item[(2)] $r=0$ and $H_0=0$: from \eqref{quelleuguaglianzesotto} it is easy to see that $\cE\geq0$ and by \eqref{1434} $\cE(M_j)=0$ for every $j\in \N\setminus\{0\}$, from which we obtain that $\cE$ is minimized on spheres. We also notice that $\cE$ is minimized on flat surfaces ($H=K=0$).
\end{enumerate}
The construction above can adapted to the constrained case by taking 
$$M_{\blu{j}}=\partial B_R\cup \partial B_{\rho/j}$$
for suitable $R,\rho>0$, where all the spherical surfaces are oriented with the outward normal, such that $\cA(M_{\blu{j}})=\cH^2(M_{\blu{j}})=a$ and $\cV(\Sigma_{M_{\blu{j}}})=v$. 
Then $\Sigma_{M_{\blu{j}}}\in\cX^{\blu{(0,1)}}_{A,c;a,v}(\Omega)$ and $\cE(M_{\blu{j}})$ has an expression similar to that in \eqref{1434}, so that the same conclusions above hold.

The case $r=0$ and $H_0 \ne 0$ is open and we do not have a counterexample at the moment.

\end{remark}
\begin{remark}[$\alpha_K=0$] In this case, the Canham--Helfrich functional $\cE$ in \eqref{100} reduces to the functional
\begin{equation}\label{100W0}
\cW_0(M)\coloneqq \alpha_H\int_M (H(p)-H_0)^2\,\de\cH^2(p),
\end{equation}
which is non-negative and is minimized by a (portion of a) sphere with mean curvature $H=H_0$. Moreover, if $H_0=0$, this further reduces to the Willmore functional $\cW$ in \eqref{100W}, which is again non-negative and minimized, for instance, on flat surfaces or on minimal surfaces. There is a vast literature on the Willmore functional both in the constrained and unconstrained case, see, \emph{e.g.}, \cite{KMS2014,NP2020,Pozzetta2021,RS2010,CS2022} in addition to those already mentioned in the Introduction. 

Here we observe that Lemma~\ref{770} provides the eigenvalue $2\alpha_H$ with multiplicity $1$ and the zero eigenvalue with algebraic multiplicity $8$. Moreover, it is necessary for the coercivity of $\cE$ that all the eigenvectors associated with the zero eigenvalue belong to $\widetilde\cX_y^\perp$ and this is not the case. Therefore, we cannot prove the coercivity in \eqref{1045} so that the direct method of the Calculus of Variations cannot be applied to show existence of minimizers. This suggests that the space of generalized Gauss graphs is not a good environment to study the Willmore functional $\cW$ of \eqref{100W}.

\end{remark}

\subsection{Regularity of minimizers}
We prove a regularity result for minimizers of $ \cE$. 
\begin{theorem}\label{thm_reg}
Let \eqref{costanti} hold and let $\Sigma \in \curv_2(\Omega)$ be a solution either of problem \eqref{1286} or of problem \eqref{1343} with $\partial\Sigma=0$, or of problems \eqref{1480}. 
Then $p_1\Sigma$ is $\cC^2$-rectifiable, that is there exists a countable family $\{S_j\}_{j\in\N}$ of surfaces of class $\cC^2$ in $\mathbb R^3$ such that 
\[
\mathcal H^2\bigg(p_1\Sigma \setminus \bigcup_{j\in\N}S_j\bigg)=0.
\]
\end{theorem}

\begin{proof}
We start by observing that, by \cite[Theorem~6.1]{D96}, since $\partial \Sigma=0$ and $|\Sigma_1|\ll |\Sigma_0|$, we get that $p_1\Sigma$ is the support of a $2$-dimensional curvature varifold (see the proof of \cite[Theorem~6.1]{D96} for the explicit construction). The regularity of $\Sigma$ is now a consequence of \cite[Theorem~1]{M13}.
\end{proof}

\begin{remark}
We point out that Theorem \ref{thm_reg} cannot be obtained using the Structure Theorem \cite[Theorem~2.10]{AST}, which asserts that if $\Sigma$ is a generalized Gauss graph then $p_1\Sigma$ is (only) $\cC^1$-rectifiable. 
\end{remark}

\bigskip

\noindent\textbf{Acknowledgements} 
The authors are members of the GNAMPA group (\emph{Gruppo Nazionale per l'Analisi Matematica, la Probabilit\`{a} e le loro Applicazioni}) of INdAM (\emph{Istituto Nazionale di Alta Matematica ``F.~Severi''}). 
AK is partially supported by the INdAM--GNAMPA 2023 Project \textit{Problemi variazionali per funzionali e operatori non-locali} (E53\-C22\-001\-930\-001), and by the Austrian Science Fund (FWF) project 10.55776/P35359.
The research of LL and MM fits within the scopes of the PRIN Project \emph{Geometric-Analytic Methods for PDEs and Applications} (2022SLTHCE) (E53D23005880006) of the Italian Ministry for University and Research.
MM acknowledges support from the MIUR grant \emph{Dipartimenti di Eccellenza 2018--2022} (E11G18000350001) of the Italian Ministry for University and Research.

\bibliography{KLM} 
\bibliographystyle{alpha}

\end{document}